\documentclass[12pt]{amsproc}

\usepackage[T2A]{fontenc} 
\usepackage[cp1251]{inputenc}
\usepackage[russian]{babel}
\usepackage{amssymb,epic,eepic}
\usepackage[dvips]{graphicx}
\numberwithin{equation}{section}

\newtheorem{thm}{Теорема}[section]

\newtheorem{rim}{Замечание}[section]
\newtheorem{dif}{Определение}[section]

\title[]{О показателях Гёльдера самоподобных функций}
\author{И.~А.~Шейпак,%
\address{Московский государственный университет
им.~М.~В.~Ломоносова, механико-математический факультет}
\email{iasheip@yandex.ru}
}
\thanks{Работа  поддержана грантом РНФ N.~17-11-01215, }

\newcommand{\Lip}{\operatorname{Lip}}

\begin{document}
\noindent УДК~517.984, 517.518.26
\begin{abstract}
Изучаются показатели Гёльдера аффинно-самоподобных функций на отрезке $[0;1]$. Получены явные формулы для показателей Гёльдера через параметры самоподобия.
\end{abstract}
\maketitle

{\small
\textbf{Ключевые слова: }\textit{Самоподобные функции, показатели Гёльдера}

\textbf{Key words: }\textit{Self-similar functions, H\"{o}lder} exponents}

\section{Введение}
Показатель Гёльдера является  одной из важнейших локальных характеристик гладкости функций. Особенно интересно для приложений знание показателей Гёльдера для непрерывных функций, производная которых равна нулю почти всюду. Второй класс непрерывных функций, представляющий интерес с этой точки зрения --- нигде недифференцируемые функции. Большинство известных функций, обладающих либо первым либо вторым свойством, фрактальны. Классическим примером функции, отвечающей первой ситуации служит функция Кантора. Стандартным представителем второго класса является функция Вейерштрасса.

Исследование  дифференциальных свойств фрактальных функций как правило, является непростой задачей.   Обычно удаётся получить результаты для частных примеров.

В частности, для функция Вейерштрасса, определяемой рядом
$$
f(x)=\sum_{k=1}^\infty b^{-k\alpha}\cos(b^k x),\quad b\in \mathbb{N}, \quad b>1,\quad 0< \alpha<1,
$$
в \cite{Hardy} было доказано, что она недифференцируема ни в одной точке.
В \cite{Zigm} (стр. 82) было показано, что эта функция удовлетворяет условию Гёльдера с показателем  $\alpha$, но не удовлетворяет условию Липшица ($\alpha=1$).

Показатель Гёльдера функции Кантора также вычислен в \cite{Zigm}, стр. 316 и равен $\log_3 2$.

Ещё один интересный класс функций основан на примере Т.\,Такаги  1903 г. (\cite{Tak}) непрерывной нигде недифференцируемой функции, задаваемой на отрезке $[0;1]$ в виде ряда
$$
T(x)=\sum_{k=0}^\infty \frac{s_0(2^k x)}{2^k},
$$
где $s_0(x)=\frac12-\left|x-\frac12\right|$.

Эту функцию можно включить в двупараметрическое семейство. Фиксируем натуральное число $n\geqslant 2$  и зададим функцию равномерно сходящимся рядом
\begin{equation}\label{eq:Tak}
T_{n,d}(x)=\sum_{k=0}^\infty d^k s_0(n^k x), \quad 0<d<1.
\end{equation}

При $n=2$, $d=\frac12$ получается классическая функция Такаги.

В 1930 г. эту функцию  переоткрыл (с другими параметрами) ван дер Варден, предложив ряд
$$
T_{10,1/10}(x)=\sum_{k=0}^\infty \frac{s_0(10^k x)}{10^k},
$$
т.е. в формуле \eqref{eq:Tak} положено $n=10$, $d=\frac{1}{10}$. Пример с $n=4$, $d=\frac14$  стал настолько классическим, что часто фигурирует в стандартных учебниках по анализу. Например, доказательство недифференцируемости  функции $T_{4,1/4}$ ни в какой точке отрезка $[0;1]$ можно найти в \cite{Fiht2}, пункт 444.  Привести полную библиографию, посвящённую этому классу функций, не представляется возможным. Укажем только обзоры \cite{AllKaw} и \cite{Lag}. Среди недавних работ упомянем \cite{Gal}.
Заметим с вышесказанным, что в математической литературе в качестве названия для этого класса функций используется и терминология "<функция ван дер Вардена"> и "<функция  Такаги"> (иногда Такаги--Ландсберга). Из-за формы графика функции $T_{2,1/2}$ иногда в литературе можно встретить название  "<функция--бланманже">.

Указанные функции относятся к частным случаям фрактальных кривых, интерес к которым особенно возобновился, после того как выяснилось, что эти кривые имеют связь с теорией всплесков и масштабирующих функций. В работах \cite{Daub2}, \cite{CavDah} были получены оценки сверху на показатель Гёльдера для фрактальных кривых в $C[0;1]$. В работе \cite{Prot2} исследовались фрактальные кривые в различных функциональных пространствах и были получены формулы на показатель гладкости в этих пространствах.

Отметим, что в этих работах результаты получены в случае так называемого "<равномерного"> разбиения отрезка $[0;1]$, на котором задавалась фрактальная кривая или функция. В работах \cite{VSh2}--\cite{GSh1} рассматривались аффинно-самоподобные функции с произвольным разбиением отрезка $[0;1]$. В работе \cite{VSh2} исследовались свойства спектральных задач с весом,  являющимися обобщённой производной самоподобной функции. В работах \cite{Sh1}, \cite{GSh1} изучались такие свойства самоподобных функций как непрерывность, монотонность, ограниченность вариации.

Целью данной работы является получение показателей Гёльдера для непрерывных аффинно-самоподобных функций при достаточно произвольном выборе параметров самоподобия, обеспечивающих лишь непрерывность функции.

Статья имеет следующую структуру. $\S 2$ содержит краткую информацию о непрерывных аффинно-самоподобных функциях.  В $\S 3$ получены формулы для показателя Гёльдера для непрерывных самоподобных функциях. В $\S 4$ рассмотрены примеры, иллюстрирующие результаты $\S 3$.

\section{Непрерывные аффинно-самоподобные функции}
Напомним конструкцию аффинно-самоподобных функций. Более подробное изложение можно найти в \cite{VSh2}, \cite{Sh1}.

Пусть фиксировано натуральное число $n>1$, и заданы положительные действительные числа
$a_1$, $a_2$, \ldots, $a_n$, удовлетворяющие условию
\[
\sum\limits_{k=1}^n a_k=1.
\]

Определим числа $\alpha_1=0$, $\alpha_{k}=\sum_{j=1}^{k-1} a_j$, $k=2,3,\ldots,n+1$. Введём булев вектор $\{e_k\}_{k=1}^n$ и определим семейство аффинных преобразований отрезка $[0,1]$ на отрезки $[\alpha_{k};\alpha_{k+1}]$:
$$
S_k(x)=a_kx+\alpha_k, \quad e_k=0; \qquad S_k(x)=-a_{k}x+\alpha_{k+1}, \quad e_k=1.
$$

Рассмотрев  произвольные (пока) наборы действительных   чисел
$\{c_k\}_{k=1}^n$, $\{d_k\}_{k=1}^n$ и $\{\beta_k\}_{k=1}^n$, где \(k=1,\ldots,n\) определим аффинное преобразование функции $f$

\begin{equation}\label{eq:samop}
[G(f)](t)=\sum\limits_{k=1}^n\left(d_k\cdot f\left(S_k^{-1}(t)\right)+c_k\cdot t+\beta_k\right)\cdot\chi_{(\alpha_k;\alpha_{k+1}]},
\end{equation}
где  через $\chi_{(a,b]}$ обозначена характеристическая функция промежутка $(a,b]$. При $k=1$ необходимо брать характеристическую функцию вида
$\chi_{[\alpha_1;\alpha_2]}=\chi_{[0;\alpha_2]}$.

Несложно проверить (см.\cite{Sh1}), что отображение $G$ является сжимающим в пространстве ограниченных функций $B[0;1]$ (и в $L_\infty[0;1]$) тогда и только тогда, когда
\begin{equation}\label{eq:szimB}
\max_{1\leqslant k\leqslant n}|d_k|<1.
\end{equation}
При этом для любых ограниченных функций $f$ и $g$ справедливо
$$
\sup_{x\in[0;1]}|[G(g)](x)-[G(f)](x)|=\max_{1\leqslant k\leqslant n}|d_k|\cdot \sup_{x\in[0;1]}|f(x)-g(x)|.
$$
Следовательно, при выполнении условия \eqref{eq:szimB} существует единственная ограниченная  функция $f$, являющаяся неподвижной точкой преобразования $G$, т.е. удовлетворяющая соотношению $G(f)=f$. График такой функции на промежутке $(\alpha_k;\alpha_{k+1}]$ $(k\ne1$, при $k=1$ график функции берётся на отрезке $[0;\alpha_2]$) есть уменьшенная и, может быть, перевёрнутая  и сдвинутая на линейную функцию копия графика на всём отрезке $[0;1]$. Т.е. с точки зрения графика самоподобие функции означает, что
\begin{equation}\label{eq:samopmatr}
\begin{pmatrix}
x\\
f(x)
\end{pmatrix}=
\bigcup_{k=1}^n
\begin{pmatrix}
a_k & 0\\
\hat c_k & d_k
\end{pmatrix}
\begin{pmatrix}
t\\
f(t)
\end{pmatrix}+
\begin{pmatrix}
\alpha_k\\
\hat\beta_k
\end{pmatrix}.
\end{equation}
Заметим, что в этой формуле используются параметры $\hat c_k$ и $\hat\beta_k$, которые, вообще говоря, отличаются от параметров $c_k$ и $\beta_k$ из формулы \eqref{eq:samop}.
Укажем связь между этим параметрами. Первая строка формулы \eqref{eq:samopmatr} означает, что при $t\in[0;1]$ следует, что $x=a_kt+\alpha_k=S_k(x)$.
Вторая строка означает, что
\begin{equation}\label{eq:samop1}
f(x)=\sum_{k=1}^n \left(d_k f(t)+\hat c_k t+\hat\beta_k\right)\cdot\chi_{(\alpha_k;\alpha_{k+1}]}=\left(d_k f\left(S_k^{-1}(x)\right)+\hat c_k S_k^{-1}(x)+\hat\beta_k\right)\cdot \chi_{(\alpha_k;\alpha_{k+1}]},
\end{equation}
откуда $c_k=\frac{\hat c_k}{a_k}$, $\beta_k=\hat\beta_k-\frac{\hat c_k\alpha_k}{a_k}$.
В зависимости от обстоятельств мы будем пользоваться либо формулой \eqref{eq:samop}, либо формулой \eqref{eq:samop1}. Формулы \eqref{eq:samop1} более удобны для компьютерного моделирования графика самоподобной. Рис.1 и рис.2  выполнены на основе формулы \eqref{eq:samop1}.

В данной работе нас будет интересовать непрерывные самоподобные функции.

В работе \cite{Sh1} были получены условия, при которых $G$ является сжимающим в $C[0;1]$ в ситуации, когда $e_k=0$ $k=1,2,\ldots,n$.
Случай произвольного булева вектора  $\{e_k\}_{k=1}^n$ был рассмотрен в \cite{GSh1}. Для дальнейшего изложения структура булева вектора $\{e_k\}_{k=1}^\infty$ не важна, поэтому далее мы будем считать, что все $e_k=0$.

Нам понадобиться следующий результат.
\begin{thm} (см. \cite{Sh1}, теорема 3.1)
Сжимающее в $B[0;1]$ преобразование $G$ с булевым вектором $e_1=e_2=\ldots=e_n=0$ задает непрерывную функцию тогда и только тогда, когда выполнены
следующие условия:
\begin{equation}\label{eq:condD}
\max_{1\leqslant k\leqslant n}|d_k|<1,
\end{equation}
\begin{gather}\label{eq:contf2_b}
\hat \beta_1=f(0)(1-d_1),\\
\hat\beta_k=\sum_{j=1}^{k-1} \hat c_j+f(1)\sum_{j=1}^{k-1}\label{eq:condBC}
d_j+f(0)(1-\sum_{j=1}^{k}d_k),\quad k=2,3,\ldots,n,\\
\sum_{j=1}^{n}\hat c_j+(f(1)-f(0))\sum_{j=1}^{n}d_j=f(1)-f(0).\label{eq:contf2_e}
\end{gather}
\end{thm}

Непрерывная неподвижная  для отображения $G$ функция $f$  называется \emph{самоподобной} в $C[0;1]$. Набор чисел $\{a_k\}$, $\{c_k\}$ (или $\{\hat c_k\}$), $\{d_k\}$, $\{\beta_k\}$ (или $\{\hat \beta_k\}$) и $\{e_k\}$ называют \emph{параметрами самоподобия} отображения $G$. В дальнейшем условимся, что если параметры  $\{c_k\}$ или $\{e_k\}$ не указаны, то они предполагаются нулевыми.

В дальнейшем мы будем рассматривать только непрерывные самоподобные функции. Не ограничивая общности, положим $f(0)=0$. Из формул \eqref{eq:condD}, \eqref{eq:contf2_b} следует, что $f(0)=0$ тогда и только тогда, когда $\hat\beta_1=0$. Что касается значения $f(1)$, то можно рассматривать два класса функций.  Первый: $f(0)\ne f(1)$. В этом случае, не ограничивая общности, можно считать, что $f(1)=1$. Второй класс: $f(1)=f(0)=0$.

\section{Непрерывность по Гёльдеру самоподобных функций.}\label{HolderExp}
Напомним определение.

\begin{dif}
Функция $f$, заданная на отрезке $[0;1]$, удовлетворяет условию Гёльдера с показателем $\alpha\in(0;1]$ если существует константа $C\geqslant 0$, такая что
\begin{equation}\label{def:Hold}
\sup_{x,y\in[0;1],x\ne y}\dfrac{|f(x)-f(y)|}{|x-y|^\alpha}\leqslant C.
\end{equation}
\end{dif}

Заметим, что если функция удовлетворяет условию Гёльдера с показателем $\alpha_0\in(0;1]$, то эта функция также удовлетворяет этому условию при любом $\alpha<\alpha_0$. Условие Гёльдера при $\alpha=1$ также называют \textit{условием Липшица}, а  функцию, удовлетворяющую  условию Гёльдера при $\alpha=1$ называют \textit{липшицевой}. Свойство \eqref{def:Hold} часто называют непрерывностью по Гёльдеру с показателем $\alpha$. Для краткости введём обозначение
$$
M(f,\alpha):=\sup_{x,y\in[0;1],x\ne y}\dfrac{|f(x)-f(y)|}{|x-y|^\alpha}.
$$

C функциями, непрерывными по Гёльдеру с показателем $\alpha$, естественным образом связывают банаховы пространства гёльдеровых функций $C^{0,\alpha}[0;1]$ и $C^{0,1}[0;1]=\Lip[0;1]$, снабженные нормами
$\|f\|_{C^{0,\alpha}}=\|f\|_{C}+M(f,\alpha)$ и $\|f\|_{C^{0,1}}=\|f\|_{C}+M(f,1)$ соответственно. Здесь $\|f\|_{C}$ --- стандартная норма в пространстве непрерывных функций $\|f\|_{C}:=\max_{x\in [0;1]}|f(x)|$. Пространство $C^{0,1}[0;1]$ называют липшицевым.

Из уравнения самоподобия \eqref{eq:samop} следует, что для всех $x\in [0;1]$ и всякого $j=1,2,\ldots n$ самоподобная функция $f$ удовлетворяет соотношению
\begin{equation}\label{eq:S}
f(S_j(x))=d_j \cdot f(x)+c_j\cdot S_j(x)+\beta_j.
\end{equation}

Из этого соотношения и определения отображений $S_j$ (при $e_j=0$) следует, что для всех $x,y\in[0;1]$ выполнено
\begin{equation}\label{eq:razn}
f(S_j(x))-f(S_j(y))=d_j\cdot\left(f(x)-f(y)\right)+c_ja_j(x-y).
\end{equation}

\begin{thm}\label{thm:d>a1}
Пусть  существует индекс $i$ такой, что $|d_i|>a_i$. Тогда для самоподобной функции $f\in C[0;1]$  выполнено  условие
Гёльдера при всех $\alpha$, удовлетворяющих условию
$$
\alpha\leqslant \min_{1\leqslant j\leqslant n}\dfrac{\ln |d_{j|}}{\ln a_{j}}.
$$
\end{thm}
\begin{proof}

Точки $\alpha_1=0$ и $\alpha_{n+1}=1$ будем называть точками нулевого уровня. Определим индуктивно множество $A_k$ точек $k$-го уровня.
Если построено множество $A_{k-1}$ точек уровня $k-1$, $A_{k-1}=\{x_i\}_{i=1}^{m_{k-1}}$, то по определению множество $A_k=\{S_j(x_i),\quad j=1,2,\ldots, n,\quad x_i\in A_{k-1}, i=1,2,\ldots,m_{k-1}\}$.

Несложно видеть, что множество $k$-го уровня ($k=0,1,\ldots$) содержит $n^{k}+1$ точку (т.е. $m_k=n^{k}+1$) и $A_{k}\subset A_{k+1}$. Расстояние между двумя соседними точками $k$-го уровня не превосходит $(\displaystyle\max_{1\leqslant j\leqslant n}a_j)^k\to 0$ при $k\to \infty$. Следовательно, множество $A:=\cup_{k=1}^\infty A_k$  всюду плотно на  отрезке $[0;1]$. Множество $A$ инвариантно относительно преобразований $S_j$, $j=1,2,\ldots,n$, т.е. $S_j(A)=A$, $j=1,2,\ldots,n$, где $S_j(A)=\left\{S_j(x_i), \quad x_i\in A\right\}$.

Непрерывная самоподобная функция  $f$ является равномерным пределом кусочно-линейных непрерывных функций $f_k$ c узлами из множества $A_k$. Для самоподобной функции $f$, удовлетворяющей условиям $f(1)=1$, $f(0)=0$, возьмём $f_0(x)=x$. Для  самоподобной функции $f$, удовлетворяющей условиям $f(1)=f(0)=0$, возьмём $f_0(x)\equiv 0$. Функции $f_k$ определяются индуктивно $f_k=G(f_{k-1})$, $k=1,2,\ldots$
Введём обозначение: $S_{j_1,j_2,\ldots j_N}(x)=S_{j_N}\circ S_{j_{N-1}}\circ\ldots S_{j_1}(x)$.

Из вышесказанного следует, что для самоподобной функции $f$
\begin{multline*}
\sup_{x,y\in[0;1],x\ne y}\dfrac{|f(x)-f(y)|}{|x-y|^\alpha}=\sup_{x,y\in A,x\ne y}\dfrac{|f(x)-f(y)|}{|x-y|^\alpha}=\\=
\limsup_{N\to\infty}\max_{1\leqslant j\leqslant N}\dfrac{|f(S_{j_1,j_2,\ldots j_N}(x))-f(S_{j_1,j_2,\ldots j_N}(y))|}{|S_{j_1,j_2,\ldots j_N}(x)-S_{j_1,j_2,\ldots j_N}(y)|^\alpha}.
\end{multline*}
При этом в правой части последнего равенства можно брать произвольные точки $x,y\in A$, $x\ne y$.

Заметим, что длина отрезка $|S_{j_1,j_2,\ldots j_N}(x)-S_{j_1,j_2,\ldots j_N}(y)|=a_{j_1}a_{j_2}\ldots a_{j_N}|x-y|$ не зависит от структуры булева вектора $\{e_k\}_{k=1}^n$.

Положим $x=1$, $y=0$.  Используя соотношение \eqref{eq:razn}, вычислим $f(S_{j_1,j_2,\ldots j_N}(x))-f(S_{j_1,j_2,\ldots j_N}(y))$:
\begin{multline*}
f(S_{j_1,j_2,\ldots j_N}(x))-f(S_{j_1,j_2,\ldots j_N}(y))=d_{j_1}d_{j_2}\ldots d_{j_N}(f(x)-f(y))+(x-y)\left(c_{j_1}d_{j_2}\ldots d_{j_N}+\right.\\+
\left.a_{j_1}c_{j_2}d_{j_3}\ldots d_{j_N}+\ldots+a_{j_1}a_{j_2}\ldots a_{j_{N-2}}c_{j_{N-1}}d_{j_N}+a_{j_1}a_{j_2}\ldots a_{j_{N-2}}a_{j_{N-1}}c_{j_N}\right).
\end{multline*}
Введя обозначения $q=\displaystyle \max_{1\leqslant j\leqslant n}\dfrac{|d_j|}{a_j^\alpha}$, $c=\displaystyle \max_{1\leqslant j\leqslant n}\dfrac{|c_j|}{a_j}$, $a=\displaystyle \max_{1\leqslant j\leqslant n} a_j$, получим
\begin{multline}\label{mult:Otz}
\dfrac{|f(S_{j_1,j_2,\ldots j_N}(x))-f(S_{j_1,j_2,\ldots j_N}(y))|}{|S_{j_1,j_2,\ldots j_N}(x)-S_{j_1,j_2,\ldots j_N}(y)|}\leqslant
2\|f\|_C\cdot q^N+\\+
c|x-y|^{1-\alpha}\cdot(q^{N-1}+q^{N-2}a^{1-\alpha}+\ldots+q(a^{1-\alpha})^{N-2}+(a^{1-\alpha})^{N-1})=\\=
2\|f\|_C\cdot q^N+c|x-y|^{1-\alpha}\cdot\dfrac{q^N-(a^{1-\alpha})^{N-1}}{q-a^{1-\alpha}}.
\end{multline}

В силу того, что $a\in(0;1)$, $\alpha\in(0;1)$, величина справа ограничена тогда и только тогда, когда
$$
q\leqslant 1\Leftrightarrow \max_{1\leqslant j\leqslant n}\dfrac{|d_j|}{a_j^\alpha}\leqslant 1 \Leftrightarrow \alpha\leqslant \min_{1\leqslant j\leqslant n}\dfrac{\ln |d_{j|}}{\ln a_{j}}.
$$

\end{proof}

Результат для показателей Гёльдера для класса функций, для которых при всех $i=1,2,\ldots,n$ выполнено $|d_i|>a_i$ получен в \cite{SingZaidl} (Следствие 5.10). Такое сильное ограничение авторам понадобилось, чтобы показать недифференцируемость ни в одной точке отрезка $[0;1]$ самоподобных функций с такими параметрами самоподобия.

\begin{thm}\label{thm:d>a2}
Пусть  существует индекс $i$ такой, что $|d_i|>a_i$. Тогда для самоподобной функции $f\in C[0;1]$  не выполнено условие Гёльдера при
$\alpha>\displaystyle\min_{1\leqslant j\leqslant n} \dfrac{\ln |d_j|}{\ln a_j}$.
\end{thm}
\begin{proof}
Обозначим через $i_0$ тот индекс, для которого $\dfrac{\ln |d_{i_0}|}{\ln a_{i_0}}=\min_{1\leqslant j\leqslant n} \dfrac{\ln |d_j|}{\ln a_j}$ (если таких индексов несколько, то для определённости в качестве $i_0$ выберем наименьший из них).

Рассмотрим два случая. 1) При всех $i=1,2,\ldots,n$ $c_i=0$. Тогда

Тогда
\begin{multline*}
M(f,\alpha)\geqslant \sup_{x,y\in[\alpha_{i_0};\alpha_{i_0+1}], x\ne y}\dfrac{|f(x)-f(y)|}{|x-y|^\alpha}=\\=
\sup_{\tilde x,\tilde y\in [0;1]\tilde x\ne \tilde y}\dfrac{|f(S_{i_0}(\tilde x))-f(S_{i_0}(\tilde y))|}{|S_{i_0}(\tilde x)-S_{i_0}(\tilde y)|^\alpha}=
\sup_{\tilde x,\tilde y\in [0;1]\tilde x\ne \tilde y}\dfrac{|d_{i_0}\cdot (f(\tilde x)-f(\tilde y))|}{a_{i_0}^\alpha|\tilde x-\tilde y|^\alpha}\geqslant
\dfrac{|d_{i_0}|}{a_{i_0}^\alpha}M(f,\alpha).
\end{multline*}

Итерируя это неравенство $N$ раз, получим
$$
M(f,\alpha)\geqslant\left(\dfrac{|d_{i_0}|}{a_{i_0}^\alpha}\right)^NM(f,\alpha),
$$
где $\left(\dfrac{|d_{i_0}|}{a_{i_0}^\alpha}\right)^N\to \infty$ при $N\to\infty$ и $\alpha>\dfrac{\ln |d_{i_0}|}{\ln a_{i_0}}$. Поскольку $|d_{i_0}|>a_{i_0}$, то самоподобная функция $f$ отлична от постоянной, поэтому $M(f,\alpha)$ неограничена.

2) Есть коэффициенты $c_i\ne 0$. Заметим, что если при некотором $i$  $c_i\ne 0$, то на отрезке $[\alpha_i;\alpha_{i+1}]$ самоподобная функция равна сумме своей уменьшенной и сдвинутой версии и линейной функции. Но линейная добавка не увеличивает показатель Гёльдера.
\end{proof}

\begin{thm}\label{thm:d=a}
Пусть при $i=1,2,\ldots,n$ для параметров самоподобия функции $f\in C[0;1]$ выполнены условия $\dfrac{|d_i|}{a_i}=1$.
Тогда самоподобная  функция $f$  удовлетворяет условию Гёльдера при всех  $\alpha\in(0;1)$.
\end{thm}
\begin{proof}
Рассуждения аналогичны, проведённым при доказательстве теоремы \ref{thm:d>a1}. В силу того, что при всех $i=1,2,\ldots,n$, $|d_i|=a_i$, $0<a_i<1$, то при всех
$\alpha\in(0;1)$ получаем, что и в этом случае $q=\displaystyle\max_{1\leqslant j\leqslant n}\dfrac{|d_j|}{a_j^\alpha}<1$. Сохраняя те же обозначения для $q$, $c$ и $a$ и используя формулу \eqref{mult:Otz}, получаем, что полунорма  $M(f,\alpha)$ ограничена при любом $\alpha\in(0;1)$.
\end{proof}

\begin{rim}
В случае, когда при всех $i=1,2,\ldots,n$ выполнены равенства $|d_i|=a_i$ и функция $f$ не является линейной, то нестрого неравенства, как в теореме \ref{thm:d>a1}, для показателя $\alpha$ получить нельзя. 
\end{rim}

\begin{thm}\label{thm:d<a}
Если при  $i=1,2,\ldots,n$ выполнены неравенства $|d_i|<a_i$, то  самоподобная функция $f$  удовлетворяет  условию
Гёльдера при всех $\alpha\leqslant 1$ (т.е. функция $f$ является липшицевой).
\end{thm}
\begin{proof}
При условиях теоремы в формуле \eqref{mult:Otz} $q=\displaystyle\max_{1\leqslant j\leqslant 1}\dfrac{d_j}{a_j^\alpha}<1$ при всех $\alpha\leqslant 1$. Следовательно полунорма $M(f,\alpha)$ ограничена при всех таких значениях $\alpha$.
\end{proof}

\section{Примеры.}

\subsection{Функция Кантора.}

Функция Кантора является самоподобной для следующего набора параметров самоподобия: $n=3$, $a_1=a_2=a_3=\frac13$, $c_1=c_2=c_3=0$, $d_1=d_3=\frac12$, $d_2=0$, $\beta_1=0$, $\beta_2=\beta_3=\frac12$. Следует уточнить, что при $d_2=0$ полагается $\ln d_2=-\infty$. В соответствии с теоремой \ref{thm:d>a1}, показатель Гёльдеровости  функции Кантора равен $\log_3 2$.

\subsection{Семейство функций Такаги--Ландсберга.}
Функцию семейства Такаги--Ландсберга мы определим при натуральном $n\geqslant 2$ и $d\in (0;1)$ как сумму ряда
$$
T_{n,d}(x)=\sum_{k=0}^\infty d^k s_0(n^k x), \quad s_0(x)=\frac12-\left|x-\frac12\right|, \quad x\in[0;1].
$$
Несложно проверить, что $T_{n,d}$ удовлетворяет функциональному уравнению
$$
T_{n,d}(x)=s_0(x)+d\cdot T_{n,d}(nx).
$$
Если $n$ чётно, то из этого соотношения и определения функции $s_0$ следует, что эта функция является самоподобной со следующими параметрами самоподобия:
\begin{gather*}
a_j=\frac1n, \quad d_j=d, \quad j=1,2,\ldots,n\\
\hat\beta_j=\frac{j-1}{n}, \quad j=1,2,\ldots, \frac{n}{2}+1,\quad  \beta_{j}=\hat\beta_{\frac{n}{2}+1}-\frac{j-1}{n}\quad j=\frac{n}{2}+2,\frac{n}{2}+3,\ldots, n,\\
\hat c_j=-\hat c_{n-j+1}=\frac{1}{n}, \quad j=1,2,\ldots, \frac{n}{2}.
\end{gather*}

Поскольку $T_{n,d}(0)=T_{n,d}(1)=0$, то  параметры $b_j$ на самом деле восстанавливаются по параметрам $c_j$ по формулам \eqref{eq:condBC}. Выбирая различные параметры $n$ и $d$, можно получить примеры для каждой из теорем \ref{thm:d>a1}, \ref{thm:d=a}, \ref{thm:d<a}.

а) Класс функций Такаги $T_{n,1/n}$ ($n\geqslant 2$, $a=d=\frac1n$) подпадает под условие теоремы \ref{thm:d=a}, поэтому удовлетворяет условию Гёльдера при любом $\alpha\in(0;1)$. C другой стороны, эта функция не является липшицевой, т.к. в соответствии с теоремой Радемахера (см., например, \cite{Fed},  (теорема 3.1.6)) липшицева функция имеет производную почти всюду, а функции Такаги не дифференцируемы при $\frac{1}{n}\leqslant d<1$ ни в одной точке отрезка $[0;1]$.

б) Функция Такаги $T_{2,1/4}$ не только лишицева, но и принадлежит пространству $C^1[0;1]$ поскольку совпадает на отрезке $[0;1]$ с функцией $g(x)=2x-2x^2$. В этом легко убедиться, проверив соотношение $G(g)=g$, в котором  для $G$ используются параметры самоподобия функции $T_{2,1/4}$. Действительно, в соответствии с формулой \eqref{eq:samop1} и параметрами самоподобия для $T_{2,1/4}$ имеем:
\begin{gather*}
[G(g)](x)=\frac14(4x-8x^2)+\frac12(2x)=2x-2x^2\quad\text{при}\quad x\in\left[0;\frac12\right],\\
[G(g)](x)=\frac14(2(2x-1)-2(2x-1)^2)-\frac12(2x-1)+\frac12=2x-2x^2\quad\text{при}\quad x\in\left[\frac12;1\right].
\end{gather*}

Пример  графика функции $T_{4,1/4}$ приведён на рис. 1.

\begin{picture}(360,190)
\put(0,30){
\begin{picture}(180,150)
\put(0,0){\includegraphics[scale=0.15]{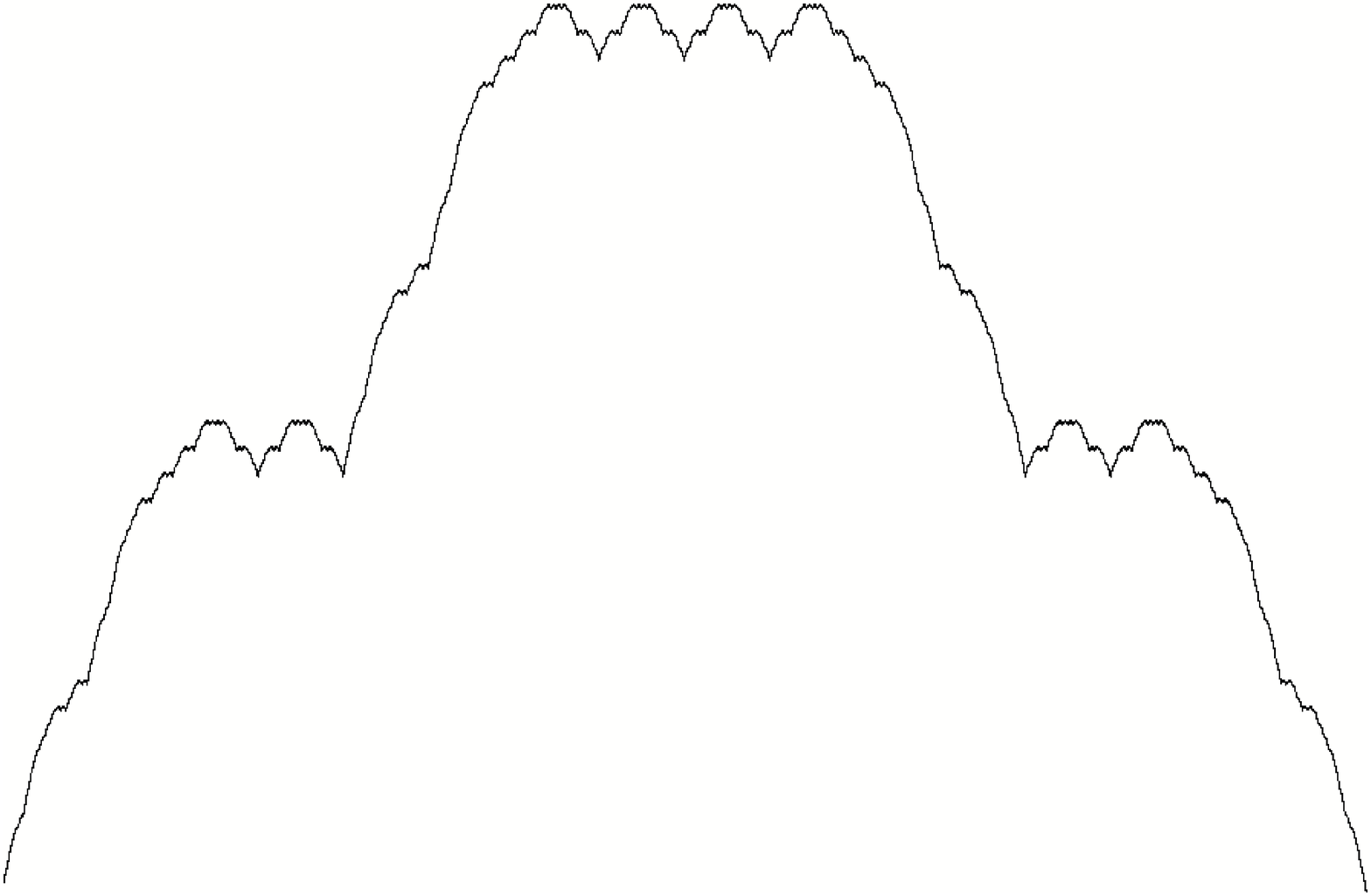}}
\put(0,0){\vector(1,0){150}}
\put(0,0){\vector(0,1){140}}
\thinlines
\multiput(134,0)(0,10){9}{\line(0,1){3}}
\multiput(0,88)(9,0){15}{\line(1,0){3}}
\put(132,-10){\small 1}
\put(-18,85){\tiny 8/15}
\put(5,-15){\small Рис. 1. Функция $T_{4,1/4}$}
\end{picture}}
\put(200,30){
\begin{picture}(180,150)
\put(-0.5,0){\includegraphics[scale=0.2]{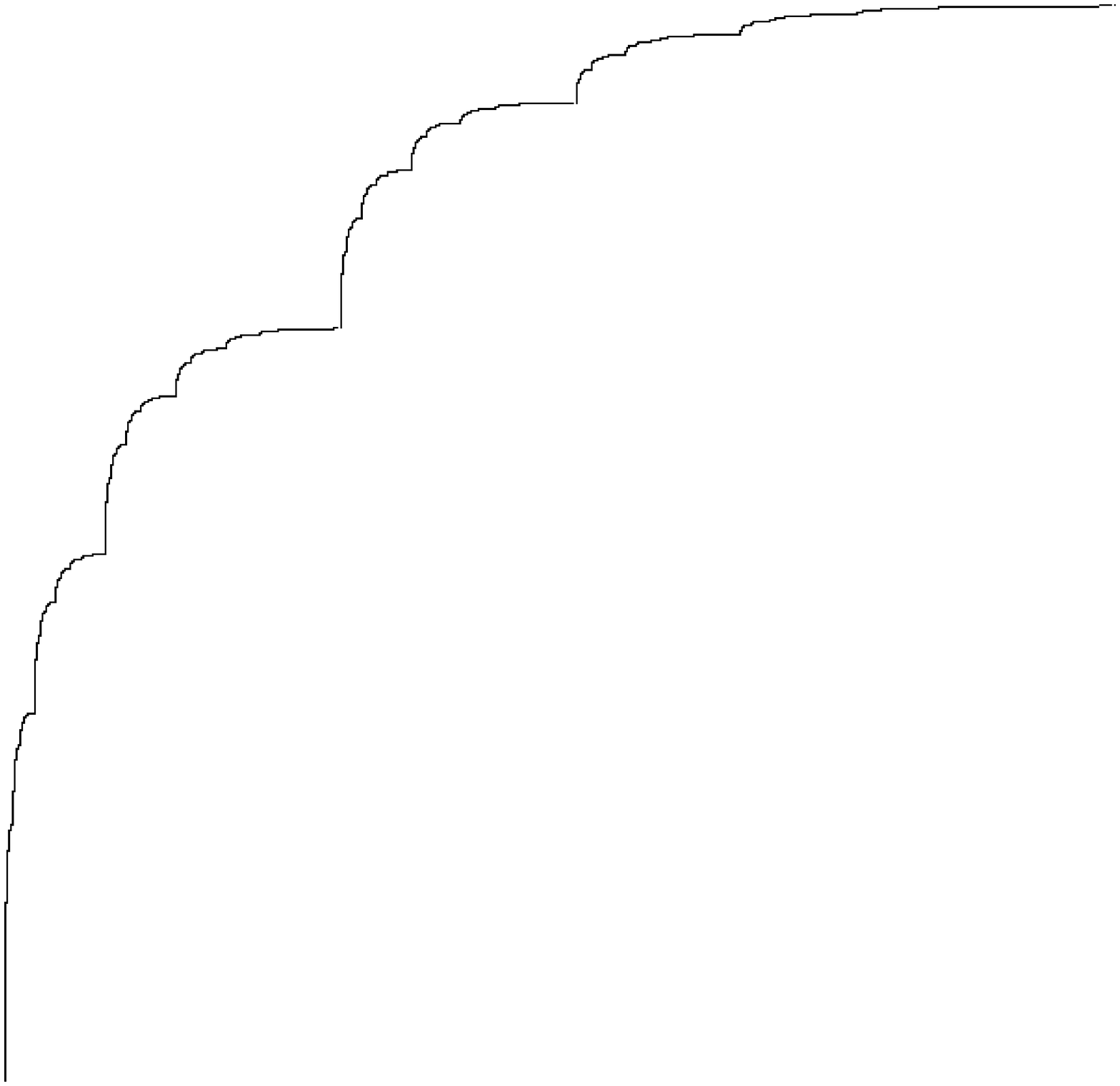}}
\put(0,0){\vector(1,0){150}}
\put(0,0){\vector(0,1){140}}
\multiput(119,0)(0,10){12}{\line(0,1){3}}
\multiput(0,118)(8.8,0){14}{\line(1,0){3}}
\put(121,-10){\small 1}
\put(-7,114){\small 1}
\put(5,-15){\small Рис. 2. Функция $\mathcal{S}_{0.3,0.7}$}
\end{picture}}
\end{picture}

\begin{picture}(360,190)
\put(0,30){
\begin{picture}(180,150)
\put(-0.5,0){\includegraphics[scale=0.2]{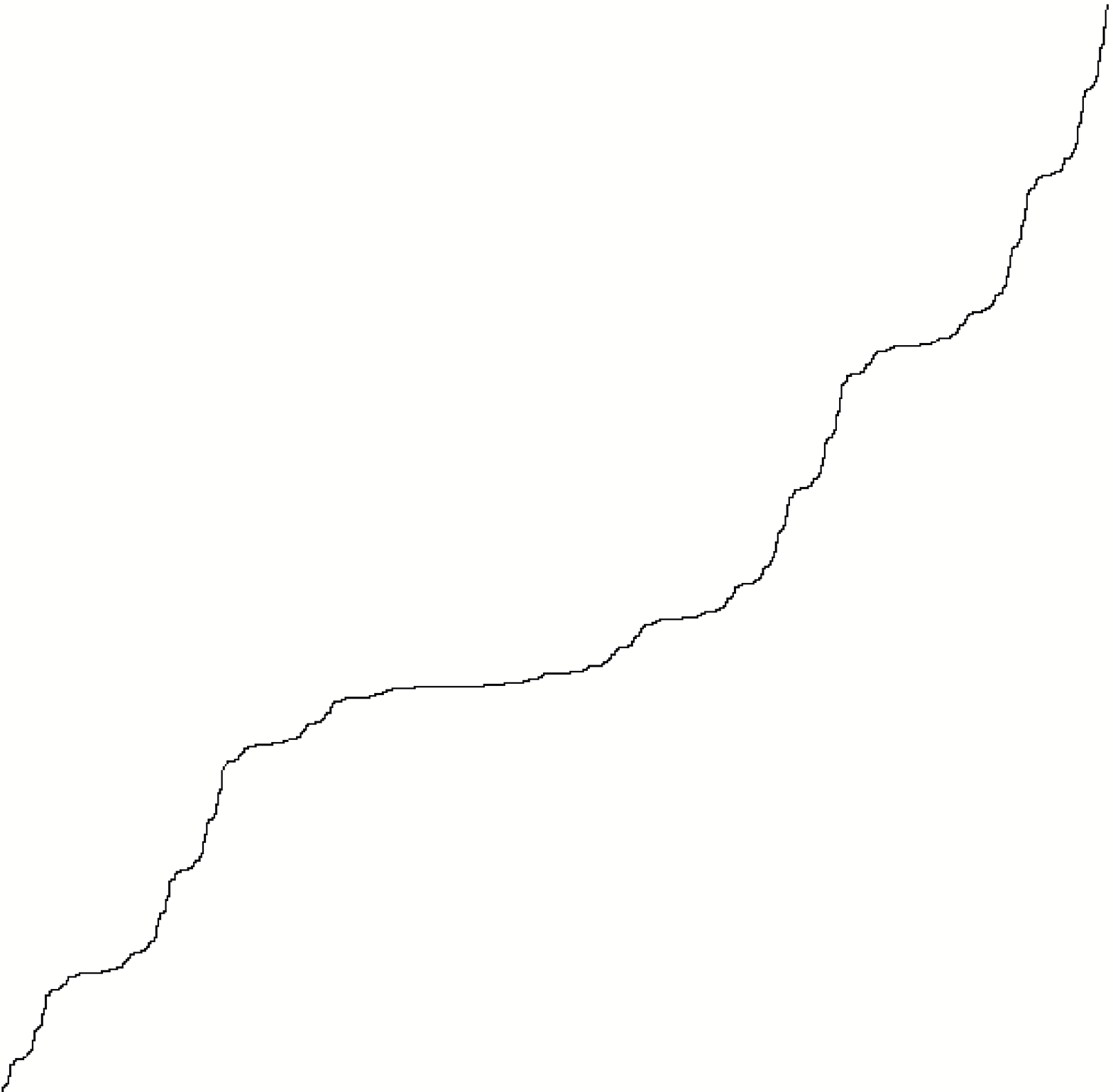}}
\put(0,0){\vector(1,0){150}}
\put(0,0){\vector(0,1){140}}
\multiput(119,0)(0,10){12}{\line(0,1){3}}
\multiput(0,118)(8.8,0){14}{\line(1,0){3}}
\put(121,-10){\small 1}
\put(-7,114){\small 1}
\put(5,-15){\small Рис. 3. Функция $\mathcal{S}_{3,\{a\},\{d\}}$}
\end{picture}}
\put(200,30){
\begin{picture}(180,150)
\put(0,0){\includegraphics[scale=0.15]{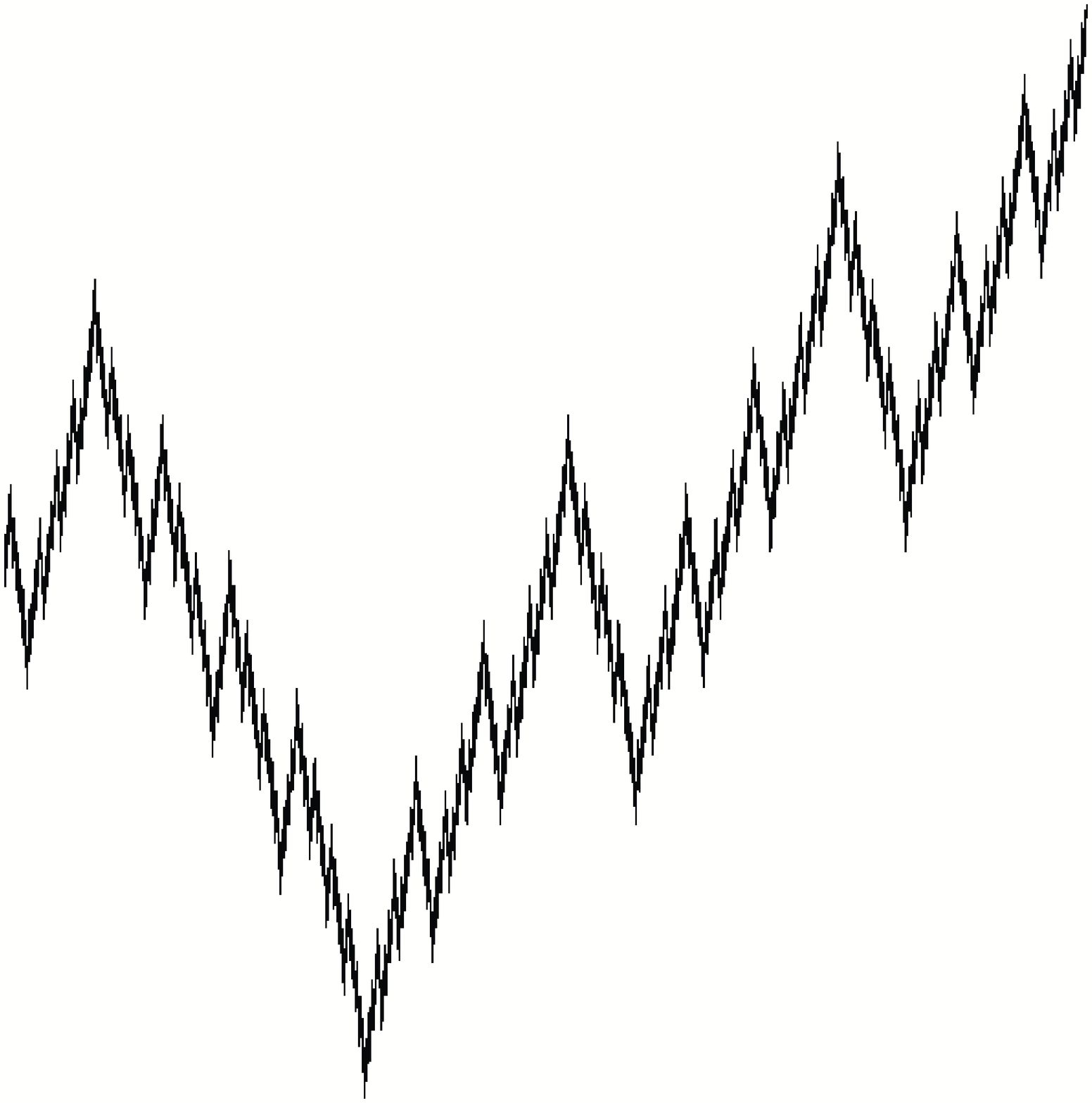}}
\put(0,35){\vector(1,0){150}}
\put(0,0){\vector(0,1){140}}
\thinlines
\multiput(90,35)(0,9){6}{\line(0,1){3}}
\multiput(0,92)(9,0){10}{\line(1,0){3}}
\multiput(0,0)(9,0){4}{\line(1,0){3}}
\put(-18,-2){\tiny -1/2}
\put(90,28){\tiny 1}
\put(-6,88){\tiny 1}
\put(15,-18){\small Рис. 4. Функция $K$}
\end{picture}}
\end{picture}

\subsection{Функции Салема.}

Обобщённая  функция Салема  $\mathcal{S}_{a,d}$ является неподвижной  для  сжимающего в $C[0;1]$ отображения   при следующих параметрах самоподобия  $0<a<1$, $0<d<1$, $d\ne \frac12$, $a\ne d$:
$$
G(f)(x)=\left(d\cdot f\left(\frac{x}{a}\right)\right)\cdot \chi_{[0;a]}+\left((1-d)\cdot f\left(\frac{x-a}{1-a}\right)+d\right)\chi_{[a;1]}.
$$
Другими словами, $n=2$, $a_1=a$, $a_2=1-a$, $\beta_1=0$, $d_1=d$, $d_2=1-d$, $\beta_2=d$, $d\ne \frac12$.
Салем в своей статье \cite{Salem} рассматривал случай "<равномерного"> разбиения отрезка $[0;1]$: $a=\frac12$. Салем исследовал эту  функцию как пример строго возрастающей функции, производная которой почти всюду равна нулю.
Поскольку функция $g(x)=\frac{\ln x}{\ln(1-x)}$ убывает, то
$$
\min\left(\frac{\ln d}{\ln a}; \frac{\ln(1-d)}{\ln(1-a)}\right)=\left\{\begin{aligned}
&\frac{\ln d}{\ln a}, \; d>a,\\
&\frac{\ln(1-d)}{\ln(1-a)}, \; d<a.
\end{aligned}\right.
$$
Таким образом, по теореме \ref{thm:d>a1} показатель Гёльдера функции Салема равен $\frac{\ln d}{\ln a}$ при $d>a$ и $\frac{\ln(1-d)}{\ln(1-a)}$ при $d<a$.
При $d=a$ функция Салема вырождается в линейную функцию $f(x)=x$ и является липшицевой (и даже класса $C^1[0;1]$). Пример графика функции Салема с параметрами самоподобия $a=3/10$, $d=7/10$ приведён на рис. 2.

Класс обобщённых функций Салема можно расширить, если делить отрезок $[0;1]$ на большее число подотрезков. Точнее, для натурального числа $n\geqslant 2$
определим следующие параметры самоподобия: $a_i,d_i\in(0;1)$, $\sum_{i=1}^n a_i=\sum_{i=1}^n d_i=1$, $\beta_1=0$, $\beta_j=\sum_{i=1}^{j-1} d_i$.
Предполагается, что для всех индексов $i=1,2,\ldots,n-1$ выполнено $\frac{d_i}{a_i}\ne\frac{d_{i+1}}{a_{i+1}}$.
Функция $\mathcal{S}_{n,\{a\},\{d\}}$ определяется как неподвижная точка сжимающего отображения $G$ с указанными параметрами самоподобия. В соответствии с формулами \eqref{eq:condD}--\eqref{eq:contf2_e} она будет непрерывной, а в соответствии с \cite{Sh1} (Теорема 4.1) --- строго возрастающей (с учётом условий $d_i>0$). Показатели Гёльдера функций этого класса определяются в соответствии с теоремой \ref{thm:d>a1}. Пример графика функции, определяемой параметрами
$n=3$, $a_1=1/5$, $a_2=1/2$, $a_3=3/10$, $d_1=3/10$, $d_2=1/5$, $d_3=1/2$ приведён на рис.3.

\subsection{Функция Киссвиттера.}
Ещё один пример непрерывной нигде недифференцируемой функции $K$ был рассмотрен в \cite{KW}. Эта функция задаётся следующими параметрами самоподобия: $n=4$ $a_1=a_2=a_3=a_4=\frac14$, $-d_1=d_2=d_3=d_4=\frac12$. C учётом того, что эта функция удовлетворяет условиям $K(0)=0$, $K(1)=1$, то параметры $\hat \beta_i$ находятся из условий \eqref{eq:condBC}. Этот пример интересен тем, что один из параметров $\{d_i\}_{i=1}^n$ отрицателен (см. рис. 4). Показатель Гёльдера этой функции равен
$\frac{\ln\frac12}{\ln\frac14}=\frac12$.

\end{document}